\documentclass[11pt]{article}
\usepackage{latexsym}
\usepackage{amsfonts}
\usepackage{amssymb}
\usepackage{mathrsfs}

\usepackage{epic}
\usepackage{eepic}
\usepackage{hyperref}

\newcommand{\be}{\begin{equation}}
\newcommand{\ee}{\end{equation}}
\newcommand{\bea}{\begin{eqnarray}}
\newcommand{\eea}{\end{eqnarray}}
\newcommand{\bean}{\begin{eqnarray*}}
\newcommand{\eean}{\end{eqnarray*}}
\newcommand{\brray}{\begin{array}}
\newcommand{\erray}{\end{array}}

\newcommand{\newsection}[1]{\setcounter{equation}{0} 
\setcounter{dfn}{0}
\section{#1}}

\newtheorem{dfn}{Definition}[section]
\newtheorem{thm}[dfn]{Theorem}
\newtheorem{lmma}[dfn]{Lemma}
\newtheorem{ppsn}[dfn]{Proposition}
\newtheorem{crlre}[dfn]{Corollary}
\newtheorem{xmpl}[dfn]{Example}
\newtheorem{rmrk}[dfn]{Remark}

\newcommand{\bdfn}{\begin{dfn}}
\newcommand{\bthm}{\begin{thm}}
\newcommand{\blmma}{\begin{lmma}}
\newcommand{\bppsn}{\begin{ppsn}}
\newcommand{\bcrlre}{\begin{crlre}}
\newcommand{\bxmpl}{\begin{xmpl}}
\newcommand{\brmrk}{\begin{rmrk}}

\newcommand{\edfn}{\end{dfn}}
\newcommand{\ethm}{\end{thm}}
\newcommand{\elmma}{\end{lmma}}
\newcommand{\eppsn}{\end{ppsn}}
\newcommand{\ecrlre}{\end{crlre}}
\newcommand{\exmpl}{\end{xmpl}}
\newcommand{\ermrk}{\end{rmrk}}

\newcommand{\bbc}{\mathbb{C}}
\newcommand{\bbz}{\mathbb{Z}}
\newcommand{\bbn}{\mathbb{N}}
\newcommand{\bbr}{\mathbb{R}}

\newcommand{\cla}{\mathcal{A}}

\newcommand{\clh}{\mathcal{H}}

\newcommand{\clk}{\mathcal{K}}
\newcommand{\cll}{\mathcal{L}}

\newcommand{\prf}{\noindent{\it Proof\/}: }
\newcommand{\ots}{\otimes}
\newcommand{\raro}{\rightarrow}

\newcommand{\lgl}{\langle}
\newcommand{\rgl}{\rangle}

\newcommand{\nn}{\nonumber}

\def \qed { \mbox{}\hfill 
$\Box$\vspace{1ex}}

\begin{document}


\author{{\sc Partha Sarathi Chakraborty}\hspace{-.2em}
  \thanks{The first author 
would like to acknowledge  support from the 
National Board of Higher Mathematics, India.}\hspace{.4em} and 
{\sc Arupkumar Pal}}
\title{Equivariant  spectral triples on 
the quantum $SU(2)$ group}
\maketitle
   \begin{abstract}
We characterize all equivariant odd  spectral triples for the quantum $SU(2)$ 
group acting on its $L_2$-space and
having a nontrivial Chern character. It is shown that the dimension of 
an equivariant spectral triple is at least three, and given any element
of the $K$-homology group of $SU_q(2)$, there is an equivariant odd 
spectral triple  of dimension 3 inducing that element.
The method employed to get equivariant spectral triples in
the quantum case is then used for classical $SU(2)$, and we prove that
for $p<4$, there does not exist any equivariant spectral triple
with nontrivial $K$-homology class and dimension $p$
acting on the $L_2$-space.
 \end{abstract}
{\bf AMS Subject Classification No.:} {\large 58}B{\large 34}, {\large 46}L{\large 87}, {\large
  19}K{\large 33}\\
{\bf Keywords.} Spectral triples, quantum group.

   \newsection{Introduction} 
Study of quantum groups originated in the
   early eighties in the work of Fadeev, Sklyanin \& Takhtajan in the
   context of quantum inverse scattering theory, and was given a more 
definitive shape by Drinfeld.  It picked up
   momentum during the mid eighties, and connections were established
   with various other areas in mathematics.  They were first studied 
   in the topological setting independently by Woronowicz (\cite{wo1}) and
 Vaksman \& Soibelman (\cite{VS}), who treated the
   $q$-deformation of the $SU(2)$ group. Woronowicz then went on to
   characterize the family of compact quantum groups and studied their
   representation theory.
   
   Noncommutative geometry was introduced around the same time by
   Alain Connes, drawing inspiration mainly from the work of Atiyah
   and Kasparov.  While their classical counterparts are very
   intimately connected to each other, there has so far been very
   little work on the connection and relationship between the two
   notions of quantum groups and noncommutative geometry.  One of the
   first questions that one would like to settle, for example, is that
   given a compact quantum group, does it admit a Dirac operator that
   is equivariant under its own (co-)action.  There has been some work
   on this theme (\cite{b-k}, \cite{gos}), but none of these resolve
   the question satisfactorily.  Present article is a modest attempt
   towards answering this.
   
   The most well-known example of a compact quantum group is the
   $q$-deformation of the $SU(2)$-group, which has been studied very
   thoroughly for real values of the deformation parameter $q$ by
   Woronowicz in \cite{wo1}.  It has a natural (co-)action on itself,
   so that it can be thaught of as an $SU_q(2)$-homogeneous space. We
   investigate geometries on this homogeneous space equivariant under
   the (co-)action of $SU_q(2)$.  The earliest work on noncommutative
   geometry on $SU_q(2)$ is probably the paper by Masuda and Watanabe
   (\cite{m-w}), but their treatment was more from the point of view
   of noncommutative topology; they do not talk about the Dirac
   operator which is of fundamental importance in Connes' theory, and
   which captures topological as well as geometric information about
   the space concerned.  More recently, Bibikov \& Kulish (\cite{b-k})
   and Goswami (\cite{gos}) made attempts to get an equivariant `Dirac'
   operator on $SU_q(2)$, but none of them could accomplish it
   satisfactorily within the framework of Connes' theory, which is
   what we plan to do in the present article.
   We restrict ourselves to odd spectral triples.  
We characterize all equivariant odd spectral triples on the $L_2$ space of
the haar state.  In particular, we show the existence of a 3-summable equivariant
spectral triple. It is also shown that an equivariant spectral triple
can not be $p$-summable for $p<3$.
Then we go on to prove 
that the associated Chern character is nontrivial by computing the
pairing between the induced Fredholm module and a generator for
$K_1(C(SU_q(2)))$, which is $\bbz$.  Computation of the pairing along
with the results of Rosenberg \& Schochet (\cite{r-s}) shows that the
associated Fredholm module is a generator of $K^1$.  One immediate
corollary is the universality of equivariant odd spectral triples in
the sense that given any odd spectral triple, there is an equivariant
one that induces the same element in $K^1$.
In the last section, we show that the spectral triples that we produce
is purely a noncommutative phenomena---for classical $SU(2)$, there 
does not exist any equivariant 3-summable Dirac operator $D$ with nontrivial
sign acting on the $L_2$-space.

Except in section~5, where we treat the classical case,
we will assume $q$ to be  \emph{a real parameter
lying in the interval $(0,1)$}.

\newsection{Preliminaries} 
To fix notation, let us give here a very
brief description of the quantum $SU(2)$ group.  The $C^*$-algebra of
continuous functions on $SU_q(2)$, to be denoted by $\cla$, 
is the $C^*$-algebra generated by two elements $\alpha$
and $\beta$ satisfying the following relations:
\bean 
\alpha^*\alpha+\beta^*\beta=I,&&\alpha\alpha^*
+q^2\beta\beta^*=I,\\
\alpha\beta-q\beta\alpha=0,&&\alpha\beta^*-q\beta^*\alpha=0,\\
\beta^*\beta&=&\beta\beta^*.  
\eean 
We will denote by $\cla_f$ the dense $*$-subalgebra of $\cla$ generated
by $\alpha$ and $\beta$. Group structure is given by the coproduct 
$\Delta:\cla\raro \cla\otimes\cla$ defined  by
\bean
\Delta(\alpha)&=&\alpha\ots\alpha-q\beta^*\ots\beta,\\
\Delta(\beta)&=&\beta\ots\alpha+\alpha^*\ots\beta.  
\eean 

For two continuous linear functionals $\rho_1$ and $\rho_2$ on $\cla$,
one defines their convolution product by: 
$\rho_1\ast\rho_2(a)=(\rho_1\otimes\rho_2)\Delta(a)$.
It is known (\cite{wo1}) that $\cla$ admits a faithful state $h$, 
called the Haar state,
that satisfies
\[
h\ast \rho (a) = h(a)\rho(I) = \rho\ast h(a)
\]
for all continuous linear functionals $\rho$ and all $a\in\cla$.
We will denote by $\clh$ the GNS space associated with this state.

The representation theory of $SU_q(2)$ is strikingly similar to its
classical counterpart.  In particular, for each $n\in\{0,\frac{1}{2},
1,\ldots\}$, there is a unique irreducible unitary representation
$t^{(n)}$ of dimension $2n+1$.  Denote by $t^{(n)}_{ij}$ the
$ij$\raisebox{.4ex}{th} entry of $t^{(n)}$. These are all elements of
$\cla_f$ and they form an orthogonal basis for $\clh$. Denote by
$e^{(n)}_{ij}$ the normalized $t^{(n)}_{ij}$'s, so that
$\{e^{(n)}_{ij}: n=0,\frac{1}{2},1,\ldots, i,j=-n,-n+1,\ldots, n\}$ is
an orthonormal basis.

\brmrk
{\rm
One has to be a little careful here, because, unlike in the classical case,
the choice of matrix entries does affect the orthogonality relations. 
Therefore one has to specify the matrix entries one is working with. 
In our case, $t^{(n)}_{ij}$'s are 
the same as in Klimyk \& Schmuedgen (page~74,~\cite{kl-sch}). 
}
\ermrk

\emph{We will use the symbol $\nu$ to denote the number $1/2$} throughout 
this article, just to make some expressions occupy less space.
Using formulas for Clebsch-Gordon coefficients, and the orthogonality 
relations (page~80--81 and equation~(57), page~115 in \cite{kl-sch}), 
one can write down the actions of 
$\alpha$, $\beta$ and $\beta^*$ on $\clh$ explicitly as follows:
\bea
\alpha: e^{(n)}_{ij}  &\mapsto& 
   a_+(n,i,j) e^{(n+\nu)}_{i-\nu ,j-\nu } 
         + a_-(n,i,j)  e^{(n-\nu )}_{i-\nu ,j-\nu },\label{alpha}\\
\beta:e^{(n)}_{ij}  &\mapsto& 
    b_+(n,i,j)  e^{(n+\nu )}_{i+\nu ,j-\nu } 
         + b_-(n,i,j)  e^{(n-\nu )}_{i+\nu ,j-\nu },\label{beta}\\
\beta^*:e^{(n)}_{ij}  &\mapsto& 
    b^+_+(n,i,j)  e^{(n+\nu )}_{i-\nu ,j+\nu } 
        + b^+_-(n,i,j)  e^{(n-\nu )}_{i-\nu ,j+\nu },\label{betastar}
\eea
where
\bea
a_+(n,i,j)  & =& \Bigl(q^{2(n+i)+2(n+j)+2}
 \frac{(1-q^{2n-2j+2})(1-q^{2n-2i+2})}{(1-q^{4n+2})(1-q^{4n+4})}
                          \Bigr)^\nu,\label{aplus}\\
a_-(n,i,j)&=&\Bigl(\frac{(1-q^{2n+2j})(1-q^{2n+2i})}
                                   {(1-q^{4n})(1-q^{4n+2})}\Bigr)^\nu,\label{aminus}\\
b_+(n,i,j)&=& - \Bigl(q^{2(n+j)}\frac{(1-q^{2n-2j+2})(1-q^{2n+2i+2})}
                              {(1-q^{4n+2})(1-q^{4n+4})}\Bigr)^\nu,\label{bplus}\\
b_-(n,i,j)&=&\Bigl(q^{2(n+i)}\frac{(1-q^{2n+2j})(1-q^{2n-2i})}{(1-q^{4n})
                                                  (1-q^{4n+2})}\Bigr)^\nu,\label{bminus}\\
b^+_+(n,i,j)&=&\Bigl(q^{2(n+i)}\frac{(1-q^{2n+2j+2})(1-q^{2n-2i+2})}
                        {(1-q^{4n+2})(1-q^{4n+4})}\Bigr)^\nu,\label{bpplus}\\
b^+_-(n,i,j)&=& - \Bigl(q^{2(n+j)}
  \frac{(1-q^{2n-2j})(1-q^{2n+2i})}{(1-q^{4n})(1-q^{4n+2})}\Bigr)^\nu.\label{bpminus}
\eea

\newsection{Equivariant spectral triples}
In this section, we will formulate the notion of equivariance, and
investigate the behaviour of $D$, where $D$ is the Dirac operator of an equivariant 
spectral triple.

In the classical context of a compact Lie group $G$, a left invariant
differential operator is one that commutes with the left regular
representation of $G$.  Now in the case of abelian $G$, the
$C^*$-algebra generated by the left regular representation is nothing
but $C(\widehat{G})$. Therefore we can rephrase the left invariance
condition as a commutation condition with $C(\widehat{G})$.  For
$C(SU_q(2))$, Woronowicz has explicitly described the generators for
$C(\widehat{G})$. Therefore, a proper analog of a left invariant Dirac
operator would be a Dirac operator commuting with these generators.

Let $A_0$ and $A_1$ be the following operators on $\clh$:
\bean
A_0\quad &:& e^{(n)}_{ij} \mapsto q^j e^{(n)}_{ij},\\
A_1\quad &:&e^{(n)}_{ij} \mapsto \cases{0 & if $j=n$,\cr
                       (q^{-2n}+q^{2n+2}-q^{-2j}-q^{2j+2})^\nu e^{(n)}_{ij+1} & if $j<n$.}
\eean
The operators $A_0$ and $A_1$ generate the $C^*$-algebra of continuous
functions on the dual of $SU_q(2)$ and thus are the `generators' of
the regular representation of $SU_q(2)$ (For more details, see
\cite{p-w}; $A_0$ and $A_1$ are the operators $\mathbf{a}$ and $\mathbf{n}$ there). We
say that an operator $T$ on $\clh$ is \textbf{equivariant} if it
commutes with $A_0$, $A_1$ and $A_1^*$. It is clear that any
equivariant self-adjoint operator with discrete spectrum must be of
the form

\be \label{dirac}
D: e^{(n)}_{ij}\mapsto d(n,i)e^{(n)}_{ij},
\ee
where $d(n,i)$'s are real.
Assume then that $D$ is such an operator. Let us first write down the
commutators of $D$ with $\alpha$ and $\beta$. 
\bea 
[D,\alpha]e^{(n)}_{ij} &=& a_+(n,i,j)(d(n+\nu,i-\nu)-d(n,i))
                                                      e^{(n+\nu)}_{i-\nu,j-\nu}\nn\\
&&{}     +a_{-}(n,i,j)(d(n-\nu,i-\nu)-d(n,i))e^{(n-\nu)}_{i-\nu,j-\nu}, \label{bdd1}\\
{}[D,\beta]e^{(n)}_{ij} &=& b_+(n,i,j)(d(n+\nu,i+\nu)-d(n,i))
                                                  e^{(n+\nu)}_{i+\nu,j-\nu}\nn\\
      &&{}    +b_{-}(n,i,j)(d(n-\nu,i+\nu)-d(n,i))e^{(n-\nu)}_{i+\nu,j-\nu}.\label{bdd2}
\eea
We  are now in a position to prove the following.
\bppsn \label{condn1}
Let $D$ be an operator of the form  $e^{(n)}_{ij}\mapsto d(n,i)e^{(n)}_{ij}$.
Then $[D,a]$ is bounded for all $a\in\cla_f$ if and only if $d(n,i)$'s satisfy
the following two conditions:
\begin{equation}\label{bnd1}
d(n+\nu,i+\nu)-d(n,i)=O(1),
\end{equation}
\be \label{bnd3}
d(n+\nu,i-\nu)-d(n,i)  =  O(n+i+1).
\ee
\eppsn
\prf
Assume that $[D,a]$ is bounded for all $a\in\cla_f$. Then,
in particular, $[D,\alpha]$ and $[D,\beta]$ are bounded, so that
there is a positive constant $C$ such that
\begin{displaymath}
\|[D,\alpha]\|\leq C,\quad  \|[D,\beta]\|\leq C.
\end{displaymath}
It follows from equations~(\ref{bdd1}) and (\ref{bdd2}) that 
\begin{equation}\label{bdd3}
|a_+(n,i,j)(d(n+\nu,i-\nu)-d(n,i))|^2 + |a_{-}(n,i,j)(d(n-\nu,i-\nu)-d(n,i))|^2\leq C^2,
\end{equation}
\begin{equation}\label{bdd4}
|b_+(n,i,j)(d(n+\nu,i+\nu)-d(n,i))|^2 + |b_{-}(n,i,j)(d(n-\nu,i+\nu)-d(n,i))|^2\leq C^2
\end{equation}
for all $n$, $i$ and $j$.
From the second inequality above, we get 
\begin{displaymath}
|b_+(n,i,j)(d(n+\nu,i+\nu)-d(n,i))| \leq C \quad \forall n, i, j.
\end{displaymath}
Now
\begin{displaymath}
|b_+(n,i,j)|=\left(\frac{q^{2n+2j}-q^{4n+2}}{1-q^{4n+2}}\right)^\nu
             \left(\frac{1-q^{2n+2i+2}}{1-q^{4n+4}}\right)^\nu.
\end{displaymath}
Hence
\begin{displaymath}
1-q^2\leq \frac{1-q^2}{1-q^{4n+4}}\leq \max_{j}|b_+(n,i,j)|^2
    =\frac{1-q^{2n+2i+2}}{1-q^{4n+4}}\leq \frac{1}{1-q^4}.
  \end{displaymath}
Hence
$|d(n+\nu,i+\nu)-d(n,i)|\leq \frac{C}{(1-q^2)^{1/2}}$ for all $n$, $i$,
i.\ e.\ we have (\ref{bnd1}).
We also have from equation~(\ref{bdd3}), 
$|a_+(n,i,j)(d(n+\nu,i-\nu)-d(n,i))|\leq C$.
But 
\begin{displaymath}
a_+(n,i,j)=q\left(\frac{q^{2n+2j}-q^{4n+2}}{1-q^{4n+2}}\right)^\nu
             \left(\frac{q^{2n+2i}-q^{4n+2}}{1-q^{4n+4}}\right)^\nu.
\end{displaymath}
Hence
\begin{displaymath}
\max_j|a_+(n,i,j)|=q\left(\frac{q^{2n+2i}-q^{4n+2}}{1-q^{4n+4}}\right)^\nu.
\end{displaymath}
Therefore 
\begin{displaymath}
q\left(\frac{q^{2n+2i}-q^{4n+2}}{1-q^{4n+4}}\right)^\nu |d(n+\nu,i-\nu)-d(n,i)|\leq C\quad
\forall n,i.
\end{displaymath}
Consequently,
$q^{n+i}|d(n+\nu,i-\nu)-d(n,i)|\leq q^{-1}C\frac{1}{(1-q^2)^\nu}$,
i.\ e.\ 
\begin{equation}\label{bnd2}
|d(n+\nu,i-\nu)-d(n,i)| =  O(q^{-n-i}).
\end{equation}
Let us next write the difference $d(n+\nu,i-\nu)-d(n,i)$
as follows:
\bean
\lefteqn{\sum_{r=0}^{n+i-1} (d(n+\nu-r\nu,i-\nu-r\nu)-
                                                  d(n+\nu-(r+1)\nu,i-\nu-(r+1)\nu))}\\
&&{} -\sum_{r=0}^{n+i-1} (d(n-r\nu,i-r\nu)-
                                                  d(n-(r+1)\nu,i-(r+1)\nu))\\
&& {}+d(n+\nu-(n+i)\nu,i-\nu-(n+i)\nu)-d(n-(n+i)\nu,i-(n+i)\nu)
\eean
Using this expression together with (\ref{bnd2}) for the case $n+i=0$ and  (\ref{bnd1}),
we get (\ref{bnd3}).

Next assume that the $d(n,i)$'s satisfy the conditions~(\ref{bnd1}) and (\ref{bnd3}).
We will show that  $[D,\alpha]$ and $[D,\beta]$ are bounded, which in turn will
ensure that $[D,a]$ is bounded for all $a\in\cla_f$.
It follows from~(\ref{bnd1}) and (\ref{bnd3}) that there is a positive constant $C>0$
such that 
\begin{displaymath}
|d(n+\nu,i+\nu)-d(n,i)|\leq C,\quad
q^{n+i}|d(n+\nu,i-\nu)-d(n,i)|\leq C.
\end{displaymath}
It follows from the above two inequalities that
\begin{eqnarray*}
|a_+(n,i,j)(d(n+\nu,i-\nu)-d(n,i))|&\leq & C(1-q^4)^{-1/2} ,\\
|a_{-}(n,i,j)(d(n-\nu,i-\nu)-d(n,i))| &\leq & Cq^{-1}(1-q^2)^{-1/2}.
\end{eqnarray*}
We now conclude from (\ref{bdd1}) that $[D,\alpha]$ is bounded.
Proof of boundedness of $[D,\beta]$ is similar.\qed

Next, we exploit the condition that $D$ must have compact resolvent.
It is straightforward to see that a necessary and sufficient condition
for an operator $D$  of the form  $e^{(n)}_{ij}\mapsto d(n,i)e^{(n)}_{ij}$
to have compact resolvent is that if we write the $d(n,i)$'s in a single sequence,
it should not have any limit point other than $\infty$ or $-\infty$.
As we shall see below, in presence of (\ref{bnd1}) and (\ref{bnd3}),
we can say much more about the $d(n,i)$'s. In particular, we can extract 
information about the sign of $D$ also.
\bppsn\label{condn2}
Let $D$ be an operator of the form  $e^{(n)}_{ij}\mapsto d(n,i)e^{(n)}_{ij}$
such that $d(n,i)$'s satisfy
conditions (\ref{bnd1}) and (\ref{bnd3})  and $D$ has compact resolvent. Then
\be\label{what}
\left.
\begin{minipage}{350pt}
\begin{enumerate}
\item For each $k\in\bbn$, there exists an $r_k\in\bbn$, $r_k\geq k$ such
  that $d(n, n-k)$ are of the same sign for all $n\geq r_k$.
 \item There exists an $r\in\bbn$ such that for all $k\geq r$ and for
  all $n$, $d(n,n-k)$ are of the same sign.
\end{enumerate}
\end{minipage}
\right\}
\ee
\eppsn 
\prf 
In the following diagram, each dot stands for a $d(n,i)$,
the dot at the $i$\raisebox{.4ex}{th} row and $j$\raisebox{.4ex}{th}
column representing $d(\frac{i+j}{2},\frac{j-i}{2})$ (here $i$ and $j$
range from 0 onwards).

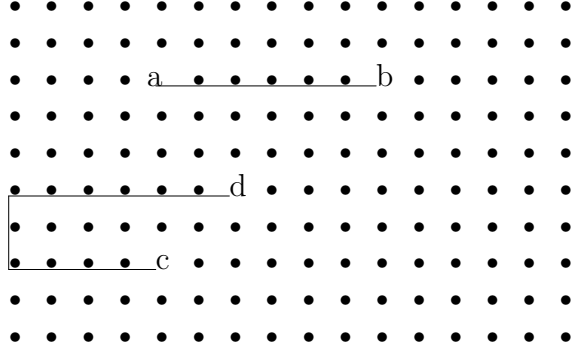
\begin{center}
\setlength{\unitlength}{0.00063333in}
\begingroup\makeatletter\ifx\SetFigFont\undefined%
\gdef\SetFigFont#1#2#3#4#5{%
  \reset@font\fontsize{#1}{#2pt}%
  \fontfamily{#3}\fontseries{#4}\fontshape{#5}%
  \selectfont}%
\fi\endgroup%
{\renewcommand{\dashlinestretch}{30}
\begin{picture}(4666,2881)(0,-10)
\put(83,2783){{\footnotesize $\bullet$}}
\put(383,2783){{\footnotesize $\bullet$}}
\put(683,2783){{\footnotesize $\bullet$}}
\put(983,2783){{\footnotesize $\bullet$}}
\put(1583,2783){{\footnotesize $\bullet$}}
\put(1883,2783){{\footnotesize $\bullet$}}
\put(2183,2783){{\footnotesize $\bullet$}}
\put(2483,2783){{\footnotesize $\bullet$}}
\put(2783,2783){{\footnotesize $\bullet$}}
\put(3083,2783){{\footnotesize $\bullet$}}
\put(3383,2783){{\footnotesize $\bullet$}}
\put(3683,2783){{\footnotesize $\bullet$}}
\put(3983,2783){{\footnotesize $\bullet$}}
\put(4283,2783){{\footnotesize $\bullet$}}
\put(4583,2783){{\footnotesize $\bullet$}}
\put(83,2483){{\footnotesize $\bullet$}}
\put(383,2483){{\footnotesize $\bullet$}}
\put(683,2483){{\footnotesize $\bullet$}}
\put(983,2483){{\footnotesize $\bullet$}}
\put(1283,2483){{\footnotesize $\bullet$}}
\put(1583,2483){{\footnotesize $\bullet$}}
\put(1883,2483){{\footnotesize $\bullet$}}
\put(2183,2483){{\footnotesize $\bullet$}}
\put(2483,2483){{\footnotesize $\bullet$}}
\put(2783,2483){{\footnotesize $\bullet$}}
\put(3083,2483){{\footnotesize $\bullet$}}
\put(3383,2483){{\footnotesize $\bullet$}}
\put(3683,2483){{\footnotesize $\bullet$}}
\put(3983,2483){{\footnotesize $\bullet$}}
\put(4283,2483){{\footnotesize $\bullet$}}
\put(4583,2483){{\footnotesize $\bullet$}}
\put(83,2183){{\footnotesize $\bullet$}}
\put(383,2183){{\footnotesize $\bullet$}}
\put(683,2183){{\footnotesize $\bullet$}}
\put(983,2183){{\footnotesize $\bullet$}}
\put(1583,2183){{\footnotesize $\bullet$}}
\put(1883,2183){{\footnotesize $\bullet$}}
\put(2183,2183){{\footnotesize $\bullet$}}
\put(2483,2183){{\footnotesize $\bullet$}}
\put(2783,2183){{\footnotesize $\bullet$}}
\put(3383,2183){{\footnotesize $\bullet$}}
\put(3683,2183){{\footnotesize $\bullet$}}
\put(3983,2183){{\footnotesize $\bullet$}}
\put(4283,2183){{\footnotesize $\bullet$}}
\put(4583,2183){{\footnotesize $\bullet$}}
\put(83,1883){{\footnotesize $\bullet$}}
\put(83,1583){{\footnotesize $\bullet$}}
\put(83,1283){{\footnotesize $\bullet$}}
\put(83,983){{\footnotesize $\bullet$}}
\put(83,683){{\footnotesize $\bullet$}}
\put(83,383){{\footnotesize $\bullet$}}
\put(83,83){{\footnotesize $\bullet$}}
\put(383,1883){{\footnotesize $\bullet$}}
\put(383,1583){{\footnotesize $\bullet$}}
\put(383,1283){{\footnotesize $\bullet$}}
\put(383,983){{\footnotesize $\bullet$}}
\put(383,683){{\footnotesize $\bullet$}}
\put(383,383){{\footnotesize $\bullet$}}
\put(383,83){{\footnotesize $\bullet$}}
\put(683,1883){{\footnotesize $\bullet$}}
\put(683,1583){{\footnotesize $\bullet$}}
\put(683,1283){{\footnotesize $\bullet$}}
\put(683,983){{\footnotesize $\bullet$}}
\put(683,683){{\footnotesize $\bullet$}}
\put(683,383){{\footnotesize $\bullet$}}
\put(683,83){{\footnotesize $\bullet$}}
\put(983,1883){{\footnotesize $\bullet$}}
\put(983,1583){{\footnotesize $\bullet$}}
\put(983,1283){{\footnotesize $\bullet$}}
\put(983,983){{\footnotesize $\bullet$}}
\put(983,683){{\footnotesize $\bullet$}}
\put(983,383){{\footnotesize $\bullet$}}
\put(983,83){{\footnotesize $\bullet$}}
\put(1283,1883){{\footnotesize $\bullet$}}
\put(1283,1583){{\footnotesize $\bullet$}}
\put(1283,1283){{\footnotesize $\bullet$}}
\put(1283,983){{\footnotesize $\bullet$}}
\put(1283,383){{\footnotesize $\bullet$}}
\put(1283,83){{\footnotesize $\bullet$}}
\put(1583,83){{\footnotesize $\bullet$}}
\put(1883,83){{\footnotesize $\bullet$}}
\put(2183,83){{\footnotesize $\bullet$}}
\put(2483,83){{\footnotesize $\bullet$}}
\put(2783,83){{\footnotesize $\bullet$}}
\put(3083,83){{\footnotesize $\bullet$}}
\put(3383,83){{\footnotesize $\bullet$}}
\put(3683,83){{\footnotesize $\bullet$}}
\put(3983,83){{\footnotesize $\bullet$}}
\put(4283,83){{\footnotesize $\bullet$}}
\put(4583,83){{\footnotesize $\bullet$}}
\put(4583,383){{\footnotesize $\bullet$}}
\put(4583,683){{\footnotesize $\bullet$}}
\put(4583,983){{\footnotesize $\bullet$}}
\put(4583,1283){{\footnotesize $\bullet$}}
\put(4583,1583){{\footnotesize $\bullet$}}
\put(4583,1883){{\footnotesize $\bullet$}}
\put(4283,1883){{\footnotesize $\bullet$}}
\put(4283,1583){{\footnotesize $\bullet$}}
\put(4283,1283){{\footnotesize $\bullet$}}
\put(4283,983){{\footnotesize $\bullet$}}
\put(4283,683){{\footnotesize $\bullet$}}
\put(4283,383){{\footnotesize $\bullet$}}
\put(3983,383){{\footnotesize $\bullet$}}
\put(3983,683){{\footnotesize $\bullet$}}
\put(3983,983){{\footnotesize $\bullet$}}
\put(3983,1283){{\footnotesize $\bullet$}}
\put(3983,1583){{\footnotesize $\bullet$}}
\put(3983,1883){{\footnotesize $\bullet$}}
\put(3683,1883){{\footnotesize $\bullet$}}
\put(3683,1583){{\footnotesize $\bullet$}}
\put(3683,1283){{\footnotesize $\bullet$}}
\put(3683,983){{\footnotesize $\bullet$}}
\put(3683,683){{\footnotesize $\bullet$}}
\put(3683,383){{\footnotesize $\bullet$}}
\put(3383,383){{\footnotesize $\bullet$}}
\put(3083,383){{\footnotesize $\bullet$}}
\put(2783,383){{\footnotesize $\bullet$}}
\put(2483,383){{\footnotesize $\bullet$}}
\put(2183,383){{\footnotesize $\bullet$}}
\put(1883,383){{\footnotesize $\bullet$}}
\put(1583,383){{\footnotesize $\bullet$}}
\put(1583,683){{\footnotesize $\bullet$}}
\put(1883,683){{\footnotesize $\bullet$}}
\put(2183,683){{\footnotesize $\bullet$}}
\put(2483,683){{\footnotesize $\bullet$}}
\put(2783,683){{\footnotesize $\bullet$}}
\put(3083,683){{\footnotesize $\bullet$}}
\put(3383,683){{\footnotesize $\bullet$}}
\put(3383,983){{\footnotesize $\bullet$}}
\put(3383,1283){{\footnotesize $\bullet$}}
\put(3383,1583){{\footnotesize $\bullet$}}
\put(3383,1883){{\footnotesize $\bullet$}}
\put(3083,1883){{\footnotesize $\bullet$}}
\put(3083,1583){{\footnotesize $\bullet$}}
\put(3083,1283){{\footnotesize $\bullet$}}
\put(3083,983){{\footnotesize $\bullet$}}
\put(2783,983){{\footnotesize $\bullet$}}
\put(2783,1283){{\footnotesize $\bullet$}}
\put(2783,1583){{\footnotesize $\bullet$}}
\put(2783,1883){{\footnotesize $\bullet$}}
\put(2483,1883){{\footnotesize $\bullet$}}
\put(2483,1583){{\footnotesize $\bullet$}}
\put(2483,1283){{\footnotesize $\bullet$}}
\put(2483,983){{\footnotesize $\bullet$}}
\put(2183,983){{\footnotesize $\bullet$}}
\put(2183,1283){{\footnotesize $\bullet$}}
\put(2183,1583){{\footnotesize $\bullet$}}
\put(2183,1883){{\footnotesize $\bullet$}}
\put(1883,1883){{\footnotesize $\bullet$}}
\put(1583,1883){{\footnotesize $\bullet$}}
\put(1583,1583){{\footnotesize $\bullet$}}
\put(1583,1283){{\footnotesize $\bullet$}}
\put(1583,983){{\footnotesize $\bullet$}}
\put(1883,983){{\footnotesize $\bullet$}}
\put(1883,1583){{\footnotesize $\bullet$}}
\put(1283,2783){{\footnotesize $\bullet$}}
\drawline(1283,2183)(3083,2183)
\drawline(1883,1283)(83,1283)(83,683)(1283,683)
\put(1208,2183){\makebox(0,0)[lb]{\smash{{{\SetFigFont{12}{14.4}{\rmdefault}{\mddefault}{\updefault}a}}}}}
\put(3083,2183){\makebox(0,0)[lb]{\smash{{{\SetFigFont{12}{14.4}{\rmdefault}{\mddefault}{\updefault}b}}}}}
\put(1283,683){\makebox(0,0)[lb]{\smash{{{\SetFigFont{12}{14.4}{\rmdefault}{\mddefault}{\updefault}c}}}}}
\put(1883,1283){\makebox(0,0)[lb]{\smash{{{\SetFigFont{12}{14.4}{\rmdefault}{\mddefault}{\updefault}d}}}}}
\end{picture}
}
\end{center}

There are two restrictions imposed on these numbers, given by
equations (\ref{bnd1}) and (\ref{bnd3}).  Equation~(\ref{bnd1}) says
that: (i)$\,$the difference of two consecutive numbers along any row is
bounded by a fixed constant, and (\ref{bnd3}) says that: (ii)$\,$the
difference of two consecutive numbers along the $j$\raisebox{.4ex}{th}
column is $O(j+1)$. Suppose $C$ is a big enough constant which works
for both (i) and (ii).

Now suppose $a$ and $b$ are two elements in the same row. Connect them
with a path as in the diagram.  If $a$ and $b$ are of opposite sign,
then because of restriction~(i) above, there has to be some dot
between $a$ and $b$ for which the corresponding $d(n,i)$ lies in
$[-C,C]$.  Therefore, if the signs of the $d(n,i)$'s change infinitely
often along a row, one can produce infinitely many $d(n,i)$'s in the
interval $[-C,C]$. But this will prevent $D$ from having a compact
resolvent. This proves part~1.

For part~2, employ a similar argument, this time connecting two dots,
say $c$ and $d$, by a path as shown in the diagram, and observing that
the difference between any two consecutive numbers along the path is
bounded by $C$.\qed

Let $m$ and $n$ be two nonnegative integers.
Let
\begin{eqnarray*}
F(m,n)&=&\left\{d\left(\frac{j+i}{2},\frac{j-i}{2}\right): 
                        0\leq i\leq m, 0\leq j\leq n\right\},\\
S(m,n,r)&=&\left\{d\left(\frac{j+r}{2},\frac{j-r}{2}\right):  
                              j> n\right\},\quad 0\leq r\leq m,\\
T(m)&=&\left\{d\left(\frac{j+i}{2},\frac{j-i}{2}\right): i>m, j\geq 0\right\}.
\end{eqnarray*}
In the following diagram, for example, $A$ is $F(2,4)$, $B$ is $T(2)$, and
$C$, $D$ and $E$ are $S(2,4,0)$, $S(2,4,1)$ and $S(2,4,2)$ respectively.
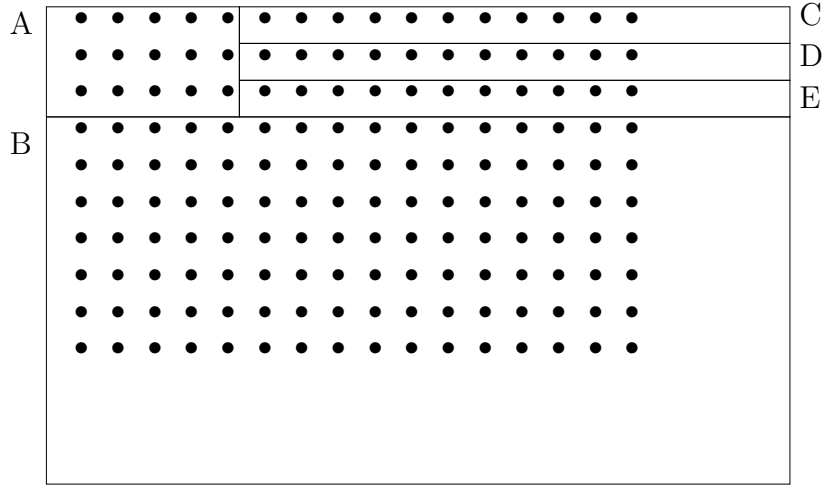
\begin{center}
\setlength{\unitlength}{0.00063333in}
\begingroup\makeatletter\ifx\SetFigFont\undefined%
\gdef\SetFigFont#1#2#3#4#5{%
  \reset@font\fontsize{#1}{#2pt}%
  \fontfamily{#3}\fontseries{#4}\fontshape{#5}%
  \selectfont}%
\fi\endgroup%
{\renewcommand{\dashlinestretch}{30}
\begin{picture}(6663,3939)(0,-10)
\put(600,3762){$\bullet$}
\put(900,3762){$\bullet$}
\put(1200,3762){$\bullet$}
\put(1500,3762){$\bullet$}
\put(2100,3762){$\bullet$}
\put(2400,3762){$\bullet$}
\put(2700,3762){$\bullet$}
\put(3000,3762){$\bullet$}
\put(3300,3762){$\bullet$}
\put(3600,3762){$\bullet$}
\put(3900,3762){$\bullet$}
\put(4200,3762){$\bullet$}
\put(4500,3762){$\bullet$}
\put(4800,3762){$\bullet$}
\put(5100,3762){$\bullet$}
\put(600,3462){$\bullet$}
\put(900,3462){$\bullet$}
\put(1200,3462){$\bullet$}
\put(1500,3462){$\bullet$}
\put(1800,3462){$\bullet$}
\put(2100,3462){$\bullet$}
\put(2400,3462){$\bullet$}
\put(2700,3462){$\bullet$}
\put(3000,3462){$\bullet$}
\put(3300,3462){$\bullet$}
\put(3600,3462){$\bullet$}
\put(3900,3462){$\bullet$}
\put(4200,3462){$\bullet$}
\put(4500,3462){$\bullet$}
\put(4800,3462){$\bullet$}
\put(5100,3462){$\bullet$}
\put(600,3162){$\bullet$}
\put(900,3162){$\bullet$}
\put(1200,3162){$\bullet$}
\put(1500,3162){$\bullet$}
\put(1800,3162){$\bullet$}
\put(2100,3162){$\bullet$}
\put(2400,3162){$\bullet$}
\put(2700,3162){$\bullet$}
\put(3000,3162){$\bullet$}
\put(3300,3162){$\bullet$}
\put(3600,3162){$\bullet$}
\put(3900,3162){$\bullet$}
\put(4200,3162){$\bullet$}
\put(4500,3162){$\bullet$}
\put(4800,3162){$\bullet$}
\put(5100,3162){$\bullet$}
\put(600,2862){$\bullet$}
\put(600,2562){$\bullet$}
\put(600,2262){$\bullet$}
\put(600,1962){$\bullet$}
\put(600,1662){$\bullet$}
\put(600,1362){$\bullet$}
\put(600,1062){$\bullet$}
\put(900,2862){$\bullet$}
\put(900,2562){$\bullet$}
\put(900,2262){$\bullet$}
\put(900,1962){$\bullet$}
\put(900,1662){$\bullet$}
\put(900,1362){$\bullet$}
\put(900,1062){$\bullet$}
\put(1200,2862){$\bullet$}
\put(1200,2562){$\bullet$}
\put(1200,2262){$\bullet$}
\put(1200,1962){$\bullet$}
\put(1200,1662){$\bullet$}
\put(1200,1362){$\bullet$}
\put(1200,1062){$\bullet$}
\put(1500,2862){$\bullet$}
\put(1500,2562){$\bullet$}
\put(1500,2262){$\bullet$}
\put(1500,1962){$\bullet$}
\put(1500,1662){$\bullet$}
\put(1500,1362){$\bullet$}
\put(1500,1062){$\bullet$}
\put(1800,2862){$\bullet$}
\put(1800,2562){$\bullet$}
\put(1800,2262){$\bullet$}
\put(1800,1962){$\bullet$}
\put(1800,1662){$\bullet$}
\put(1800,1362){$\bullet$}
\put(1800,1062){$\bullet$}
\put(2100,1062){$\bullet$}
\put(2400,1062){$\bullet$}
\put(2700,1062){$\bullet$}
\put(3000,1062){$\bullet$}
\put(3300,1062){$\bullet$}
\put(3600,1062){$\bullet$}
\put(3900,1062){$\bullet$}
\put(4200,1062){$\bullet$}
\put(4500,1062){$\bullet$}
\put(4800,1062){$\bullet$}
\put(5100,1062){$\bullet$}
\put(5100,1362){$\bullet$}
\put(5100,1662){$\bullet$}
\put(5100,1962){$\bullet$}
\put(5100,2262){$\bullet$}
\put(5100,2562){$\bullet$}
\put(5100,2862){$\bullet$}
\put(4800,2862){$\bullet$}
\put(4800,2562){$\bullet$}
\put(4800,2262){$\bullet$}
\put(4800,1962){$\bullet$}
\put(4800,1662){$\bullet$}
\put(4800,1362){$\bullet$}
\put(4500,1362){$\bullet$}
\put(4500,1662){$\bullet$}
\put(4500,1962){$\bullet$}
\put(4500,2262){$\bullet$}
\put(4500,2562){$\bullet$}
\put(4500,2862){$\bullet$}
\put(4200,2862){$\bullet$}
\put(4200,2562){$\bullet$}
\put(4200,2262){$\bullet$}
\put(4200,1962){$\bullet$}
\put(4200,1662){$\bullet$}
\put(4200,1362){$\bullet$}
\put(3900,1362){$\bullet$}
\put(3600,1362){$\bullet$}
\put(3300,1362){$\bullet$}
\put(3000,1362){$\bullet$}
\put(2700,1362){$\bullet$}
\put(2400,1362){$\bullet$}
\put(2100,1362){$\bullet$}
\put(2100,1662){$\bullet$}
\put(2400,1662){$\bullet$}
\put(2700,1662){$\bullet$}
\put(3000,1662){$\bullet$}
\put(3300,1662){$\bullet$}
\put(3600,1662){$\bullet$}
\put(3900,1662){$\bullet$}
\put(3900,1962){$\bullet$}
\put(3900,2262){$\bullet$}
\put(3900,2562){$\bullet$}
\put(3900,2862){$\bullet$}
\put(3600,2862){$\bullet$}
\put(3600,2562){$\bullet$}
\put(3600,2262){$\bullet$}
\put(3600,1962){$\bullet$}
\put(3300,1962){$\bullet$}
\put(3300,2262){$\bullet$}
\put(3300,2562){$\bullet$}
\put(3300,2862){$\bullet$}
\put(3000,2862){$\bullet$}
\put(3000,2562){$\bullet$}
\put(3000,2262){$\bullet$}
\put(3000,1962){$\bullet$}
\put(2700,1962){$\bullet$}
\put(2700,2262){$\bullet$}
\put(2700,2562){$\bullet$}
\put(2700,2862){$\bullet$}
\put(2400,2862){$\bullet$}
\put(2100,2862){$\bullet$}
\put(2100,2562){$\bullet$}
\put(2100,2262){$\bullet$}
\put(2100,1962){$\bullet$}
\put(2400,1962){$\bullet$}
\put(2400,2262){$\bullet$}
\put(2400,2562){$\bullet$}
\put(1800,3762){$\bullet$}
\drawline(375,3912)(1950,3912)(1950,3012)
        (375,3012)(375,3912)
\drawline(1950,3912)(6450,3912)(6450,3612)
        (1950,3612)(1950,3912)
\drawline(1950,3612)(6450,3612)(6450,3312)
        (1950,3312)(1950,3612)
\drawline(1950,3312)(6450,3312)(6450,3012)
        (1950,3012)(1950,3312)
\drawline(375,3012)(6450,3012)(6450,12)
        (375,12)(375,3012)
\put(75,3687){\makebox(0,0)[lb]{\smash{{{\SetFigFont{12}{14.4}{\familydefault}{\mddefault}{\updefault}A}}}}}
\put(75,2712){\makebox(0,0)[lb]{\smash{{{\SetFigFont{12}{14.4}{\familydefault}{\mddefault}{\updefault}B}}}}}
\put(6525,3762){\makebox(0,0)[lb]{\smash{{{\SetFigFont{12}{14.4}{\familydefault}{\mddefault}{\updefault}C}}}}}
\put(6525,3425){\makebox(0,0)[lb]{\smash{{{\SetFigFont{12}{14.4}{\familydefault}{\mddefault}{\updefault}D}}}}}
\put(6525,3087){\makebox(0,0)[lb]{\smash{{{\SetFigFont{12}{14.4}{\familydefault}{\mddefault}{\updefault}E}}}}}
\end{picture}
}
\end{center}
What the last proposition says is the following. There exist big enough 
integers $m$ and $n$ such that in each of the sets $T(m)$, $S(m,n,0),\ldots,S(m,n,m)$,
all elements are of the same sign, i.\ e.\ each of the sets 
$T(m)$, $S(m,n,0),\ldots,S(m,n,m)$
is contained in either $\bbr_+$ or $-\bbr_+$.
\brmrk \label{rmrk:condn4}
{\rm
One can extend the argument in the proof of the last proposition a 
little further and prove that if $D$ is as in the previous proposition, then
}
\be\label{condn4}
\left.
\begin{minipage}{350pt}
given any nonnegative real $N$, there exist positive integers $m$ and $n$ 
such that each of the sets $T(m)$, $S(m,n,0),\ldots,S(m,n,m)$
is contained in either $\{x\in\bbr: x>N\}$ or $\{x\in\bbr:x<-N\}$.
\end{minipage}
\right\}
\ee
\ermrk

\bthm
An operator $D$ on $L_2(h)$ gives rise to an equivariant spectral triple
if and only if it is of the form $e^{(n)}_{ij}\mapsto d(n,i)e^{(n)}_{ij}$,
where $d(n,i)$'s are real and satisfy conditions~(\ref{bnd1}), (\ref{bnd3}) 
and (\ref{condn4}).
\ethm
\prf
It is enough to prove that if the $d(n,i)$'s obey condition~(\ref{condn4}),
then $D$ has compact resolvent. But this is clear, because (\ref{condn4})
implies that for any real number $N>0$, the interval
$[-N,N]$ contains only a finite number of the $d(n,i)$'s.\qed

It is clear then that up to a compact perturbation, 
$D$ will have nontrivial sign if and only if the following condition holds:
\be\label{sign}
\left.
\begin{minipage}{350pt}
\begin{enumerate}
\item there exist positive integers $m$ and $n$ such that in each of the sets
  $T(m)$, $S(m,n,0),\ldots,S(m,n,m)$, all elements are of the same
  sign, and 
\item there are two sets in this collection whose elements are
  of opposite sign.
\end{enumerate}
\end{minipage}
\right\}
\ee

A natural question to ask now is whether there does indeed exist any $D$ with 
nontrivial sign satisfying (\ref{bnd1}) and (\ref{bnd3}). 
It is easy to see that the operator $D$ determined by the 
family $d(n,i)$, where
\be \label{genericd}
d(n,i)=\cases{2n+1 & if $n\neq i$,\cr
                      -(2n+1) & if $n=i$,}
\ee
satisfy all the requirements in propositions \ref{condn1} and \ref{condn2}. 
In fact, one can easily see that $D^{-3}\in\cll^{(1,\infty)}$, where 
$\cll^{(1,\infty)}$ stands for the ideal of Dixmier traceable operators. 
Thus we have the following.

\bthm
$SU_q(2)$ admits an equivariant odd  3-summable spectral triple.
\ethm

The classical $SU(2)$ has  (both topological as well as metric) dimension 3.
For $SU_q(2)$, however, the topological dimension turns out to be 1, 
as can be seen from the following short exact sequence
\begin{displaymath}
0\longrightarrow \clk\otimes C(S^1) 
\longrightarrow  \cla
\longrightarrow C(S^1) \longrightarrow 0,
\end{displaymath}
where $\clk$ denotes the algebra of compact operators.
The next theorem tells us that as far as metric dimension is concerned,
it behaves more like its classical counterpart; in fact along with the 
previous theorem, it says that the metric dimension of $SU_q(2)$ is 3.
\bthm
Let $(\cla,\clh,D)$ be an equivariant odd spectral triple. Then
$D$ can not be $p$-summable for $p<3$.
\ethm
\prf
Conditions (\ref{bnd1}) and (\ref{bnd3})
impose the following growth restriction on
the $d(n,i)$'s:
\be \label{bnd4}
\max_i |d(n,i)| = O(n).
\ee
The conclusion of the theorem follows easily from this.\qed

The next proposition gives a nice property of the operator $D$, 
namely, it says that the 
derivative of  any nonconstant function is nonzero.

\bppsn
Let $D$ be given by (\ref{genericd}). Then for $a\in\cla_f$,
$[D,a]=0$ if and only if $a$ is a scalar.
\eppsn
\prf
Take $a=\sum_{(i,j,k)\in F}c_{ijk}\alpha_i\beta^j{\beta^*}^k$, where
$F$ is a finite subset of $\bbz\times\bbn\times\bbn$ and all the $c_{ijk}$'s are 
nonzero (here $\alpha_i$ is $\alpha^i$ for $i\geq 0$ and $(\alpha^*)^{-i}$ for $i<0$). 
We will show that $[D,a]\neq 0$.

Let $m=\max \{|i|+j+k : (i,j,k)\in F\}$. Let $(r,s,t)$ be a point of 
$F$ such that $|r|+s+t=m$. Write $p=\frac{1}{2}(s-t-r)$,
$p'=\frac{1}{2}(t-s-r)$. Then it is easy to see that
\bean
\lefteqn{\lgl e^{(n+m/2)}_{pp'}, [D,a]e^{(n)}_{00}\rgl}\\
  &=& \lgl e^{(n+m/2)}_{pp'}, [D, c_{rst}\alpha_{r}\beta^{s}{\beta^*}^{t}]e^{(n)}_{00}\rgl\\
  &=& c_{rst}\prod_{i=1}^{t} b^+_+\left(n+\frac{i-1}{2}, -\frac{i-1}{2} ,\frac{i-1}{2}\right)
  \prod_{i=t+1}^{t+s}b_+\left(n+\frac{i-1}{2}, -t+\frac{i-1}{2} ,t-\frac{i-1}{2}\right)\\
&&\mbox{}\times  
  \prod_{i=s+t+1}^{m}a^{\#}_+\left(n+\frac{i-1}{2}, 
              p+\mbox{sign}(r)\frac{m-i}{2} ,p'+\mbox{sign}(r)\frac{m-i}{2}\right)\\
&&\hspace{6em}{}\times   \Bigl(d(n+m/2,p)-d(n,0)\Bigr),
\eean
where $a^{\#}_+$ stands for $a_+$ or $a^+_+$ depending on the sign of $r$.
The right hand side above is clearly nonzero because of our choice of $D$.\qed

The above proposition says, in particular, that the Dirac operator given 
by~(\ref{genericd}) is really a Dirac operator for the full tangent bundle rather than 
that of some lower dimensional subbundle.
\newsection{Nontriviality of the Chern character}

In this section we will examine the $D$ given by the family
(\ref{genericd}) in more detail and see that the nontriviality in sign
does indeed result in nontriviality at the Fredholm level. For this,
we will compute the pairing between $\mbox{sign}\,D$ and a generator
of $K_1(\cla)$. Let $u$ denote the element
$\chi_{\{1\}}(\beta^*\beta)(\beta-I)+I$ of $\cla$, where, 
for a normal operator $T$, $\chi_F(T)$ denotes the spectral
projection of $T$ corresponding to a subset $F$ of the spectrum.  
It can be shown that
this is a generator of $K_1(\cla)$.  What we will do is the following.
We will choose an invertible element $\gamma$ in $\cla_f$ that is
close enough to $u$ so that $\gamma$ and $u$ are the same in
$K_1(\cla)$.  We then compute the pairing between
$\mbox{sign}\,D$ and this $\gamma$.  
\bthm \label{main} The Chern
character of the spectral triple $(\cla_f,\clh, D)$ is nontrivial.
\ethm

\brmrk\rm
Goswami (\cite{gos}) gives an example of an equivariant operator
$D$ acting on $\clh\otimes\bbc^2$, and having bounded commutators 
with the algebra elements. But this $D$ ($|D|$ in his notation) 
is positive, hence has trivial
pairing with $K$-theory. Now 
in most cases it is possible to find a self-adjoint operator
with compact resolvent and bounded commutators with the algebra elements
just by looking at the elements affiliated to the commutant
of the algebra represented on a Hilbert space.
It is to avoid this kind of trivialities that
the nontrivial pairing is a very crucial requirement.
In the present case, as the theorem above shows,
the Dirac operator defined by (\ref{genericd})
does have a nontrivial pairing with $K$-theory.
\ermrk

Before we begin the proof of the theorem, 
let us observe from equations~(\ref{beta}) 
and~(\ref{betastar}) that
the action of $\beta\beta^*$ on $\clh$ is given by
\be\label{gmma1}
\pagebreak
(\beta\beta^*) (e^{(n)}_{ij}) = \sum_{\epsilon=-1}^1 
                       k_\epsilon(n,i,j) e^{(n+\epsilon)}_{ij}, 
\ee
where
\bea
k_1(n,i,j)  & =& -\Bigl(q^{4n+2i+2j+2}\frac{1-q^{2n+2j+2}}{1-q^{4n+2}}
                               \frac{1-q^{2n-2i+2}}{1-q^{4n+4}} \frac{1-q^{2n-2j+2}}{1-q^{4n+4}}
                       \frac{1-q^{2n+2i+2}}{1-q^{4n+6}}\Bigr)^\nu,\cr
   && \label{kplus}\\
k_0(n,i,j)  & =& q^{2(n+j)}\frac{(1-q^{2n-2j})(1-q^{2n+2i})}{(1-q^{4n})(1-q^{4n+2})}
           +q^{2(n+i)}\frac{(1-q^{2n+2j+2})(1-q^{2n-2i+2})}{(1-q^{4n+2})(1-q^{4n+4})},\cr
   && \label{kzero}\\
k_{-1}(n,i,j)  & =& -\Bigl(q^{4n+2i+2j-2}\frac{(1-q^{2n-2j})(1-q^{2n+2i})
           (1-q^{2n+2j})(1-q^{2n-2i})}{(1-q^{4n-2})(1-q^{4n})(1-q^{4n})(1-q^{4n+2})}\Bigr)^\nu.\cr
   && \label{kminus}
\eea

\noindent \textbf{Proof of theorem~\ref{main}\ }:
Choose $r\in\bbn$ such that $q^{2r}<\frac{1}{2} < q^{2r-2}$.
Define $\gamma_r=(\beta^*\beta)^r(\beta-I)+I$.
By our choice of $r$, we have
\bean
\|\gamma_r-u\| &\leq& \|(\beta^*\beta)^r-\chi_{\{1\}}(\beta^*\beta)\|\cdot\|\beta-I\|\\
                     & \leq & 2q^{2r}\;  < \;1.
\eean
Hence $\gamma_r$ and $u$ are the same in $K_1(\cla)$.
Therefore it is enough for our purpose
if we can show that the pairing between $\mbox{sign}\,D$ and $\gamma_r$ is nontrivial.
Denote by $P_k$ the projection onto the space spanned by 
$\{e^{(n)}_{n-k,j}: n,j\}$. Then $\mbox{sign}\,D=I-2P_0$.
Therefore we now want to compute the index of the operator 
$P_0\gamma_r P_0$ thaught of as an operator on $P_0\clh$.

It follows from (\ref{gmma1}) that
\be \label{gmma2}
(\beta\beta^*)^r(e^{(n)}_{ij}) = \sum_{\epsilon_t\in\{-1,0,1\}}
  \Bigl(\prod_{t=1}^r k_{\epsilon_t}(n+\sum_{s=1}^{t-1}\epsilon_s ,i,j)\Bigr)
                            e^{(n+\sum_1^r\epsilon_s)}_{ij}.
\ee
Since $\beta$ is normal, we have 
\bea
\gamma_r e^{(n)}_{ij} &=& \sum_{\epsilon_t\in\{-1,0,1\}}
  \Bigl(\prod_{t=1}^r k_{\epsilon_t}(n+\sum_{s=1}^{t-1}\epsilon_s ,i,j)\Bigr)
  \Bigl(b_+(n+\sum_1^r\epsilon_s,i,j)  e^{(n+\sum_1^r\epsilon_s+\nu )}_{i+\nu ,j-\nu } \nn\\
&&  \hspace{14em} {}+  b_-(n+\sum_1^r\epsilon_s,i,j)  
            e^{(n+\sum_1^r\epsilon_s-\nu )}_{i+\nu ,j-\nu }\Bigr)\nn\\
&&-\sum_{\epsilon_t\in\{-1,0,1\}}
  \Bigl(\prod_{t=1}^r k_{\epsilon_t}(n+\sum_{s=1}^{t-1}\epsilon_s ,i,j)\Bigr)
e^{(n+\sum_1^r\epsilon_s)}_{ij}
+e^{(n)}_{ij}. \label{gmma3}
\eea

Consequently,
\bean
\gamma_r e^{(n)}_{nj} &=& \sum_{\epsilon_t\in\{-1,0,1\}}
     \Bigl(\prod_{t=1}^r k_{\epsilon_t}(n+\sum_{s=1}^{t-1}\epsilon_s ,n,j)\Bigr)
 \Bigl(b_+(n+\sum_1^r\epsilon_s,n,j)  e^{(n+\sum_1^r\epsilon_s+\nu )}_{n+\nu ,j-\nu } \nn\\
&&  \hspace{14em} {}+  
  b_-(n+\sum_1^r\epsilon_s,n,j)  e^{(n+\sum_1^r\epsilon_s-\nu )}_{n+\nu ,j-\nu }\Bigr)\nn\\
&&-\sum_{\epsilon_t\in\{-1,0,1\}}
         \Bigl(\prod_{t=1}^r k_{\epsilon_t}(n+\sum_{s=1}^{t-1}\epsilon_s ,n,j)\Bigr)
e^{(n+\sum_1^r\epsilon_s)}_{nj} +e^{(n)}_{nj}. 
\eean

When we cut this off by $P_0$, we get 
\bean 
P_0\gamma_r e^{(n)}_{nj}
&=& \sum_{\sum\epsilon_t=0} \Bigl(\prod_{t=1}^r
k_{\epsilon_t}(n+\sum_{s=1}^{t-1}\epsilon_s ,n,j)\Bigr)
b_+(n,n,j)  e^{(n+\nu )}_{n+\nu ,j-\nu } \nn\\
&& {}+ \sum_{\sum\epsilon_t=1} \Bigl(\prod_{t=1}^r
k_{\epsilon_t}(n+\sum_{s=1}^{t-1}\epsilon_s ,n,j)\Bigr)
b_-(n+1,n,j)  e^{(n+\nu )}_{n+\nu ,j-\nu } \nn\\
&&{}-\sum_{\sum\epsilon_t=0} \Bigl(\prod_{t=1}^r
k_{\epsilon_t}(n+\sum_{s=1}^{t-1}\epsilon_s ,n,j)\Bigr)
e^{(n)}_{nj} +e^{(n)}_{nj}.  
\eean 
A closer look at the
quantities $k_\epsilon$ and $b_\pm$ tells us that if we do the
calculations modulo compact operators, which we can
because we want to compute the index, we find that there is no
contribution from the second term, while in the case of the first and
the third term, contribution comes from only the coefficient where the
product $\prod_{t=1}^r k_{\epsilon_t}(n+\epsilon_1+\ldots
+\epsilon_{t-1},n,j)$ consists solely of $k_0$'s, i.\ e.\ when each
$\epsilon_t=0$. A further examination of the terms $k_0$ and $b_+$
then yield the following:
\bean
P_0\gamma_r P_0 e^{(n)}_{nj}  &=& k_0(n,n,j)^r b_+(n,n,j) e^{(n+\nu )}_{n+\nu ,j-\nu }
                                 + (1-k_0(n,n,j)^r)e^{(n)}_{nj}\\
&=& -q^{2rn+2rj}(1-q^{2n-2j})^r q^{n+j}(1-q^{2n-2j+2})^{1/2}e^{(n+\nu )}_{n+\nu ,j-\nu }\\
      &&\hspace{4em}  + \left(1-q^{2rn+2rj}(1-q^{2n-2j})^r\right)e^{(n)}_{nj},
\eean
and
\bean
P_0\gamma_r^* P_0 e^{(n)}_{nj}  &=& 
    -q^{2rn+2rj}(1-q^{2n-2j-2})^r q^{n+j}(1-q^{2n-2j})^{1/2}e^{(n-\nu )}_{n-\nu ,j+\nu }\\
      &&\hspace{4em}      + \left(1-q^{2rn+2rj}(1-q^{2n-2j})^r\right)e^{(n)}_{nj}
\eean

From these, one can easily show that the index of $P_0\gamma_r P_0$ is $-1$.
Since $P_0$ is the eigenspace corresponding to the eigenvalue 
$-1$ of $\mbox{sign}\,D$, the value of the $K$-homology--$K$-theory pairing 
$\langle [u],[(\cla,\clh,D)] \rangle$ coming from Kasparov product of $K^1$ 
and $K_1$ is $-\mbox{index}\, P_0 \gamma_r P_0$,  which is  nonzero.\qed

\brmrk
{\rm
Strictly speaking, it is not essential to introduce the element $u$ as a 
generator for $K_1(\cla)$. It is enough if one computes the pairing between 
$\mbox{sign}\,D$ and a suitable $\gamma_r$ and show that it is nontrivial.
But the introduction of $u$ makes the choice of $\gamma_r$'s  and hence 
the proof above more transparent.
}
\ermrk

It follows from proposition~\ref{condn2} that for the purposes of
computing the index pairing, sign of any equivariant $D$ must be of
the form $I-2P$ where $P=\sum_{k\in F} P_k$, $F$ being a finite subset
of $\bbn$ ( the actual $P$ would be a compact perturbation of this). 
Conversely, given a $P$ of this form, it is easy to produce a  $D$
satisfying the conditions in proposition~\ref{condn2} for which $\mbox{sign}\,D=I-2P$.
One could, for example, take the $D$ given by $d(n,i)$'s, where
\be
d(n,i)=\cases{-(2n+1) & if $n-i \in F$,\cr
                      2n+1 & otherwise.}
\ee

We are now in a position to prove the following.

\bppsn \label{existence}
Given any $m\in\bbz$, there exists an equivariant spectral triple 
$D$ acting on $\clh$ such that 
$\langle \gamma_r,[(\cla,\clh,D)]\rangle =m$, where 
$\langle,\cdot,\cdot\rangle: K_1(\cla) \times K^1(\cla)\raro \bbz$
denotes the map coming from the Kasparov product.
\eppsn
\prf
It is enough to prove the statement for  $m$  positive.
Let $D$ be an equivariant Dirac operator whose sign is
$I-2P$ where $P=\sum_{k\in F} P_k$, $F$ being a  subset
of size $m$ of $\bbn$.
In order to compute the pairing $\langle \gamma_r, [(\cla,\clh,D)]\rangle$, we must first
have a look at $P_{k+l}\gamma_r P_k$.

We get from equation~(\ref{gmma3})
\bean
\gamma_r e^{(n)}_{n-k,j} &=& \sum_{\epsilon_t\in\{-1,0,1\}}
 \Bigl(\prod_{t=1}^r k_{\epsilon_t}(n+\sum_{s=1}^{t-1}\epsilon_s ,n-k,j)\Bigr)
         \nn\\
&&  \hspace{-3em}{}\times \Bigl(b_+(n+\sum_1^r\epsilon_s,n-k,j)  
          e^{(n+\sum_1^r\epsilon_s+\nu )}_{n-k+\nu ,j-\nu } +  
  b_-(n+\sum_1^r\epsilon_s,n-k,j)  e^{(n+\sum_1^r\epsilon_s-\nu )}_{n-k+\nu ,j-\nu }\Bigr)\nn\\
&&-\sum_{\epsilon_t\in\{-1,0,1\}}
    \Bigl(\prod_{t=1}^r k_{\epsilon_t}(n+\sum_{s=1}^{t-1}\epsilon_s ,n-k,j)\Bigr)
e^{(n+\sum_1^r\epsilon_s)}_{n-k,j} +e^{(n)}_{n-k,j}. 
\eean
and consequently,
\bean 
P_{k+l}\gamma_r e^{(n)}_{n-k,j}
&=& \sum_{\sum\epsilon_t=l} \Bigl(\prod_{t=1}^r
k_{\epsilon_t}(n+\sum_{s=1}^{t-1}\epsilon_s,n-k,j)\Bigr)
b_+(n+l,n-k,j)  e^{(n+l+\nu )}_{n-k+\nu ,j-\nu } \nn\\
&&\hspace{-4em} {}+ \sum_{\sum\epsilon_t=l+1} \Bigl(\prod_{t=1}^r
k_{\epsilon_t}(n+\sum_{s=1}^{t-1}\epsilon_s,n-k,j)\Bigr)
b_-(n+l+1,n-k,j)  e^{(n+l+\nu )}_{n-k+\nu ,j-\nu } \nn\\
&&{}-\sum_{\sum\epsilon_t=l} \Bigl(\prod_{t=1}^r
k_{\epsilon_t}(n+\sum_{s=1}^{t-1}\epsilon_s,n-k,j)\Bigr)
e^{(n+l)}_{n-k,j} +\delta_{l0}e^{(n)}_{n-k,j}.  
\eean 
Now because of the nature of the quantities $k_\epsilon$ and $b_\pm$,
 we see that for index calculations, none of the terms contribute
anything for $l\neq 0$, while for $l=0$,
the first, third and the fourth term survive, with
coefficient in the first term being
$k_0(n,n-k,j)^r b_+(n,n-k,j)$ 
and that in the third being
$(1-k_0(n,n-k,j)^r)$.
It follows from here that
\[
\mbox{index}\, P_{k}\gamma_r P_k = -1,
\]
and $P_{k+l}\gamma_r P_k$ is compact for $l\neq 0$.
Therefore the pairing between $\mbox{sign}\,D$ and $\gamma_r$ 
produces $m$.\qed

An immediate corollary of the above proposition and theorem~1.17 in \cite{r-s}
is the following universality property of equivariant spectral triples.
\bcrlre
Given any odd spectral triple $(\cla,\clk,D)$, there is an equivariant
triple $(\cla,\clh,D')$ inducing the same element in $K^1(\cla)$.
\ecrlre


Finally, we have the following characterization 
theorem for equivariant Dirac operators.
\bthm
$(\cla,\clh,D)$ is an equivariant odd spectral triple with
nontrivial Chern character if and only if $D$ is given by (\ref{dirac})
and the $d(n,i)$'s obey conditions~(\ref{bnd1}), (\ref{bnd3}),
(\ref{condn4}) and~(\ref{sign}).
\ethm
\prf
If $D$ is of the form $e^{(n)}_{ij}\mapsto d(n,i)e^{(n)}_{ij}$,
where $d(n,i)$'s are real and satisfy conditions~(\ref{bnd1}), (\ref{bnd3}),  
(\ref{condn4}) and~(\ref{sign}), then proposition~\ref{condn1} says $[D,a]$ is bounded 
and nontriviality of Chern character follows from arguments of 
proposition~\ref{existence}. Conversely,
if $D$ is equivariant, then by  
propositions~\ref{condn1}, \ref{condn2} and remark~\ref{rmrk:condn4}, we have 
(\ref{bnd1}), (\ref{bnd3}) and (\ref{condn4}).
Since $D$ has nontrivial Chern character, it has nontrivial sign 
so that we have (\ref{sign}).\qed


\newsection{The case $q=1$}
It would be interesting at this point to see what happens
in the case $q=1$, i.\ e.\ for the classical $SU(2)$.
In particular, if the operator $D$ given by (\ref{genericd})
yields anything in that case.

The representation of $C(SU(2))$ on $L_2(SU(2))$ is given by
\bean
\alpha: e^{(n)}_{ij}  &\mapsto& a_+(n,i,j) e^{(n+\nu)}_{i-\nu ,j-\nu } 
       + a_-(n,i,j)  e^{(n-\nu )}_{i-\nu ,j-\nu },\label{calpha}\\
\beta:e^{(n)}_{ij}  &\mapsto& b_+(n,i,j)  e^{(n+\nu )}_{i+\nu ,j-\nu } 
       + b_-(n,i,j)  e^{(n-\nu )}_{i+\nu ,j-\nu },\label{cbeta}
\eean
where
\bea
a_+(n,i,j)  & =& \Bigl(\frac{(n-j+1)(n-i+1)}{(2n+1)(2n+2)}
                          \Bigr)^\nu,\label{caplus}\\
a_-(n,i,j)&=&\Bigl(\frac{(n+j)(n+i)}{2n(2n+1)}\Bigr)^\nu,\label{caminus}\\
b_+(n,i,j)&=& - \Bigl(\frac{(n-j+1)(n+i+1)}{(2n+1)(2n+2)}\Bigr)^\nu,\label{cbplus}\\
b_-(n,i,j)&=&\Bigl(\frac{(n+j)(n-i)}{2n(2n+1)}\Bigr)^\nu,\label{cbminus}
\eea
\brmrk\label{qd}\rm
From these, it is immediate that the operator $D$ given by (\ref{genericd}) in this 
case does not have bounded commutators.
One should, however, note that  the 
associated Kasparov module
has a nontrivial $K$-homology class,
as can be seen by directly computing the pairing of 
$\mbox{sign}\,D$  with the fundamental unitary 
$\left(\brray{lr}\alpha &-\beta^*\cr \beta&\alpha^*\erray\right)$.
\ermrk

Observe that the representation of (the complexification of) $\mathfrak{su}(2)$ 
on $L_2(SU(2))$ is given by
\bean
h e^{(n)}_{ij} &=& (n-2j)e^{(n)}_{ij},\\
e e^{(n)}_{ij} &=& j(n-2j+1)e^{(n)}_{i,j-1},\\
f e^{(n)}_{ij} &=& e^{(n)}_{i,j+1},
\eean
where $h$, $e$ and $f$ obey
\[
[h,e]=2e,\quad [h,f]=-2f,\quad [e,f]=h.
\]
Therefore, in this case also any equivariant self-adjoint operator with discrete
spectrum must be of the form
\be\label{genericeq}
D: e^{(n)}_{ij}\mapsto d(n,i)e^{(n)}_{ij}.
\ee
Commutators of this operator with  $\alpha$ and $\beta$
are once again given by (\ref{bdd1}) and (\ref{bdd2}),
where $a_{\pm}$ and $b_{\pm}$ are now given by equations
(\ref{caplus})--(\ref{cbminus}).

\blmma\label{sign-sum}
Suppose $D$ is an operator on $L_2(SU(2))$ given by (\ref{genericeq})
and having bounded commutators with $\alpha$ and $\beta$.
Assume that except for finitely many $n$'s,
the set $\{d(n,i): i=-n,-n+1,\ldots,n\}$ contains elements of both signs,
then $D$ can not be $p$-summable for $p<4$.
\elmma
\prf
Conditions for boundedness of the commutators give us
\bea
|d(n+\nu,i+\nu)-d(n,i)| &=& O\Bigl(\left(\frac{2n+2}{n+i+1}\right)^\nu\Bigr),\label{cbdd1}\\
|d(n+\nu,i-\nu)-d(n,i)| &=& O\Bigl(\left(\frac{2n+2}{n-i+1}\right)^\nu\Bigr).\label{cbdd2}
\eea
Clearly, then, there is a $K>0$ such that
\bean
|d(n+\nu,i+\nu)-d(n,i)| &\leq& K\sqrt{n},\\
|d(n+\nu,i-\nu)-d(n,i)| &\leq& K\sqrt{n}.
\eean
By assumption, there is an $i$ such that $d(n,i)$ and $d(n,i+1)$ are of opposite signs.
Therefore, either the pair $d(n-\nu,i+\nu)$ and $d(n,i)$
or the pair $d(n-\nu,i+\nu)$ and $d(n,i+1)$ must be of opposite sign,
so that their difference is really the sum of their absolute
values. Because of the above inequalities, either $|d(n,i)|\leq K\sqrt{n}$
or $|d(n,i+1)|\leq K\sqrt{n}$.
Thus for all but finitely many $n$'s, there is an $i$ such that
$|d(n,i)|\leq K\sqrt{n}$.
Now it is a routine exercise in operator theory to show 
that $D$ can not be $p$-summable for $p<4$.\qed

\bppsn\label{sign-sum2}
Suppose $D$ be an operator on $L_2(SU(2))$ given by (\ref{genericeq}),
having nontrivial sign
and having bounded commutators with $\alpha$ and $\beta$.
Then except for finitely many $n$'s,
the set $\{d(n,i): i=-n,-n+1,\ldots,n\}$ contains elements of both signs.
\eppsn
Note that we call $\mbox{sign}\,D$ trivial if it is $I$ or $-I$ up to a compact
perturbation.

\prf
Observe from (\ref{cbdd1}) and (\ref{cbdd2}) that if we restrict ourselves to 
the region $i\geq 0$, then
\be\label{cbdd3}
|d(n+\nu,i+\nu)-d(n,i)| = O(1),
\ee
and if we restrict to $i\leq 0$, then 
\be\label{cbdd4}
|d(n+\nu,i-\nu)-d(n,i)| = O(1).
\ee
Also, it is not too difficult to see that
\be\label{cbdd5}
|d(n+1,0)-d(n,0)| = O(1).
\ee
Suppose $C>0$ is a constant that works for (\ref{cbdd1})--(\ref{cbdd5}).

For the region $i\geq 0$, using arguments similar to that employed 
in the proof of proposition~\ref{condn2} but connecting two elements
$c$ and $d$ lying on two different rows by a path as shown in the diagram here
and using (\ref{cbdd3}) and (\ref{cbdd5}),
one can show that
\begin{center}
\setlength{\unitlength}{0.00070833in}
\begingroup\makeatletter\ifx\SetFigFont\undefined%
\gdef\SetFigFont#1#2#3#4#5{%
  \reset@font\fontsize{#1}{#2pt}%
  \fontfamily{#3}\fontseries{#4}\fontshape{#5}%
  \selectfont}%
\fi\endgroup%
{\renewcommand{\dashlinestretch}{30}
\begin{picture}(4666,2881)(0,-10)
\put(83,2783){$\bullet$}
\put(383,2783){$\bullet$}
\put(683,2783){$\bullet$}
\put(983,2783){$\bullet$}
\put(1583,2783){$\bullet$}
\put(1883,2783){$\bullet$}
\put(2183,2783){$\bullet$}
\put(2483,2783){$\bullet$}
\put(2783,2783){$\bullet$}
\put(3083,2783){$\bullet$}
\put(3383,2783){$\bullet$}
\put(3683,2783){$\bullet$}
\put(3983,2783){$\bullet$}
\put(4283,2783){$\bullet$}
\put(4583,2783){$\bullet$}
\put(83,2483){$\bullet$}
\put(383,2483){$\bullet$}
\put(683,2483){$\bullet$}
\put(983,2483){$\bullet$}
\put(1283,2483){$\bullet$}
\put(1583,2483){$\bullet$}
\put(1883,2483){$\bullet$}
\put(2183,2483){$\bullet$}
\put(2483,2483){$\bullet$}
\put(2783,2483){$\bullet$}
\put(3083,2483){$\bullet$}
\put(3383,2483){$\bullet$}
\put(3683,2483){$\bullet$}
\put(3983,2483){$\bullet$}
\put(4283,2483){$\bullet$}
\put(4583,2483){$\bullet$}
\put(83,2183){$\bullet$}
\put(383,2183){$\bullet$}
\put(683,2183){$\bullet$}
\put(983,2183){$\bullet$}
\put(1583,2183){$\bullet$}
\put(1883,2183){$\bullet$}
\put(2183,2183){$\bullet$}
\put(2483,2183){$\bullet$}
\put(2783,2183){$\bullet$}
\put(3383,2183){$\bullet$}
\put(3683,2183){$\bullet$}
\put(3983,2183){$\bullet$}
\put(4283,2183){$\bullet$}
\put(4583,2183){$\bullet$}
\put(83,1883){$\bullet$}
\put(83,1583){$\bullet$}
\put(83,683){$\bullet$}
\put(83,383){$\bullet$}
\put(83,83){$\bullet$}
\put(383,1883){$\bullet$}
\put(383,1583){$\bullet$}
\put(383,1283){$\bullet$}
\put(383,983){$\bullet$}
\put(383,683){$\bullet$}
\put(383,383){$\bullet$}
\put(383,83){$\bullet$}
\put(683,1883){$\bullet$}
\put(683,1583){$\bullet$}
\put(683,1283){$\bullet$}
\put(683,983){$\bullet$}
\put(683,683){$\bullet$}
\put(683,383){$\bullet$}
\put(683,83){$\bullet$}
\put(983,1883){$\bullet$}
\put(983,1583){$\bullet$}
\put(983,1283){$\bullet$}
\put(983,983){$\bullet$}
\put(983,383){$\bullet$}
\put(983,83){$\bullet$}
\put(1283,1883){$\bullet$}
\put(1283,1583){$\bullet$}
\put(1283,983){$\bullet$}
\put(1283,383){$\bullet$}
\put(1283,83){$\bullet$}
\put(1583,83){$\bullet$}
\put(1883,83){$\bullet$}
\put(2183,83){$\bullet$}
\put(2483,83){$\bullet$}
\put(2783,83){$\bullet$}
\put(3083,83){$\bullet$}
\put(3383,83){$\bullet$}
\put(3683,83){$\bullet$}
\put(3983,83){$\bullet$}
\put(4283,83){$\bullet$}
\put(4583,83){$\bullet$}
\put(4583,383){$\bullet$}
\put(4583,683){$\bullet$}
\put(4583,983){$\bullet$}
\put(4583,1283){$\bullet$}
\put(4583,1583){$\bullet$}
\put(4583,1883){$\bullet$}
\put(4283,1883){$\bullet$}
\put(4283,1583){$\bullet$}
\put(4283,1283){$\bullet$}
\put(4283,983){$\bullet$}
\put(4283,683){$\bullet$}
\put(4283,383){$\bullet$}
\put(3983,383){$\bullet$}
\put(3983,683){$\bullet$}
\put(3983,983){$\bullet$}
\put(3983,1283){$\bullet$}
\put(3983,1583){$\bullet$}
\put(3983,1883){$\bullet$}
\put(3683,1883){$\bullet$}
\put(3683,1583){$\bullet$}
\put(3683,1283){$\bullet$}
\put(3683,983){$\bullet$}
\put(3683,683){$\bullet$}
\put(3683,383){$\bullet$}
\put(3383,383){$\bullet$}
\put(3083,383){$\bullet$}
\put(2783,383){$\bullet$}
\put(2483,383){$\bullet$}
\put(2183,383){$\bullet$}
\put(1883,383){$\bullet$}
\put(1583,383){$\bullet$}
\put(1583,683){$\bullet$}
\put(1883,683){$\bullet$}
\put(2183,683){$\bullet$}
\put(2483,683){$\bullet$}
\put(2783,683){$\bullet$}
\put(3083,683){$\bullet$}
\put(3383,983){$\bullet$}
\put(3383,1283){$\bullet$}
\put(3383,1883){$\bullet$}
\put(3083,1883){$\bullet$}
\put(3083,1583){$\bullet$}
\put(3083,1283){$\bullet$}
\put(3083,983){$\bullet$}
\put(2783,983){$\bullet$}
\put(2783,1283){$\bullet$}
\put(2783,1583){$\bullet$}
\put(2783,1883){$\bullet$}
\put(2483,1883){$\bullet$}
\put(2483,1583){$\bullet$}
\put(2483,1283){$\bullet$}
\put(2483,983){$\bullet$}
\put(2183,983){$\bullet$}
\put(2183,1283){$\bullet$}
\put(2183,1583){$\bullet$}
\put(2183,1883){$\bullet$}
\put(1883,1883){$\bullet$}
\put(1583,1883){$\bullet$}
\put(1583,1583){$\bullet$}
\put(1583,1283){$\bullet$}
\put(1583,983){$\bullet$}
\put(1883,983){$\bullet$}
\put(1883,1583){$\bullet$}
\put(1283,2783){$\bullet$}
\put(1283,1283){$\bullet$}
\put(983,683){$\bullet$}
\put(83,1283){$\bullet$}
\put(83,983){$\bullet$}
\drawline(1283,2183)(3083,2183)
\drawline(3383,1583)(1360,1583)(2260,683)(3383,683)
\put(1208,2183){\makebox(0,0)[lb]{\smash{{{\SetFigFont{10}{12.0}{\rmdefault}{\mddefault}{\updefault}a}}}}}
\put(3083,2183){\makebox(0,0)[lb]{\smash{{{\SetFigFont{10}{12.0}{\rmdefault}{\mddefault}{\updefault}b}}}}}
\put(1283,683){$\bullet$}
\put(1883,1283){$\bullet$}
\put(3383,1508){\makebox(0,0)[lb]{\smash{{{\SetFigFont{10}{12.0}{\rmdefault}{\mddefault}{\updefault}c}}}}}
\put(3383,683){\makebox(0,0)[lb]{\smash{{{\SetFigFont{10}{12.0}{\rmdefault}{\mddefault}{\updefault}d}}}}}
\end{picture}
}
\end{center}
\begin{enumerate}
\item for any given row, signs of all the $d(n,i)$'s are eventually the same,
\item there exists an integer $K>0$ such that $(K+1)$\raisebox{.4ex}{th} row onwards,
all the $d(n,i)$'s have the same sign in the region $i\geq 0$.
\end{enumerate}

Similar reasoning tells us that for the region $i\leq 0$,
\begin{enumerate}\setcounter{enumi}{2}
\item for any given column, signs of all the $d(n,i)$'s are eventually the same,
\item there exists an integer $K'>0$ such that $(K'+1)$\raisebox{.4ex}{th} column onwards,
all the $d(n,i)$'s have the same sign in the region $i\leq 0$.
\end{enumerate}
Since the intersection of the regions $i\leq0$ and $i\geq 0$ is nonempty,
it follows that 
\begin{enumerate}\setcounter{enumi}{4}
\item if we leave out the first $K$ rows and first $K'$ columns,
all the remaining $d(n,i)$'s have the same sign. 
\end{enumerate}

Let $R_{mk}$ be the set of $d(n,i)$'s 
in the $k$\raisebox{.4ex}{th} row lying on the $(m+1)$\raisebox{.4ex}{th} column onwards,
$C_{mk}$ be the set of $d(n,i)$'s in the $k$\raisebox{.4ex}{th} column
lying on the $(m+1)$\raisebox{.4ex}{th} row onwards and $T_m$ be the set of $d(n,i)$'s
at the $ij$\raisebox{.4ex}{th} position, where $i\geq m+1$ and $j\geq m+1$.
In other words, let
\bean
R_{mk} &=& \Bigl\{d(n,n-k): n\geq \frac{m+k}{2}\Bigr\},\\
C_{mk} &=& \Bigl\{d(n,k-n): n\geq \frac{m+k}{2}\Bigr\},\\
T_m &=& \Bigl\{d(n,i): n\geq m, -n+m\leq i \leq n-m\Bigr\}.
\eean 
From the observations~(1--5) above, we conclude that 
there is a big enough integer $m$ such that the sets
$R_{m0},\ldots ,R_{mm}$, $C_{m0},\ldots,C_{mm}$, and $T_m$
are all contained in either $\bbr_+$ or $-\bbr_+$,
and there are at least two sets in this collection whose elements
are of opposite signs.
This immediately tells us that for all $n\geq m$,
the set $\{d(n,i): i=-n,-n+1,\ldots,n\}$ has both 
positive as well as negative elements.\qed

Combining lemma~\ref{sign-sum} and proposition~\ref{sign-sum2}, we now get the 
following theorem.
\bthm
Suppose $(C(SU(2)),L_2(SU(2)),D)$ is an equivariant spectral triple,
and assume that $D$ has nontrivial sign.
Then $D$ can not be $p$-summable for $p<4$.
\ethm

The following example illustrates that the bound obtained in the above theorem
on the summability of $D$ is the best possible.

Let $D$ be given by the following $d(n,i)$'s 
\be
d(n,i)=\cases{ -[\sqrt{2n}] & if $i=n$,\cr
               2[\sqrt{2n}] & if $i=n-1$,\cr
               3[\sqrt{2n}] & if $i=n-2$,\cr
               \ldots &\ldots \cr
               k_{2n}[\sqrt{2n}] & if $i=n-k_{2n}+1$,\cr
               2n      & if $i\leq n-k_{2n}$,}
\ee
where $k_{2n}+1$ is the least integer greater than or equal to $\frac{2n}{[\sqrt{2n}]}$.
Then for the operator $|D|$, following are the eigenvalues along with 
their multiplicities on the $n$\raisebox{.4ex}{th} chunk,
i.\ e.\ on $\mbox{span}\{e^{(n)}_{ij}: i,j=-n,-n+1,\ldots, n\}$:
\begin{center}
\begin{tabular}{cc}
   eigenvalue  & multiplicity\\
    \hline
     $[\sqrt{2n}]$ & $2n+1$ \\
    $2[\sqrt{2n}]$ & $2n+1$ \\
     $3[\sqrt{2n}]$ & $2n+1$ \\
     \ldots &\ldots\\
     $k_{2n}[\sqrt{2n}]$ & $2n+1$ \\
     $2n$ & $(2n+1)(2n-k_{2n}+1)$
\end{tabular}
\end{center}

Therefore each integer $n\in\bbn$ is an eigenvalue for $|D|$, and the multiplicity
$m_n$ of $n$ is given by
\be
m_n=(n+1)(n+1-k_n)+\sum_{r|n}(\frac{n^2}{r^2}+1)+
       \sum_{r|n}(\frac{n^2}{r^2}+2)+\ldots +\sum_{r|n}(\frac{n^2}{r^2}+2n+1).
\ee
It follows from this that $m_n=O(n^3)$, so that $D$ is at most 4-summable.
But $D$ has nontrivial sign, and therefore by the theorem above,
it can not be $p$-summable for $p<4$. Hence $D$ is 4-summable.

Since the sign of this operator coincides with the sign of the operator given by
(\ref{genericd}), it follows from remark~\ref{qd} that
the $K$-homology class of this $D$ is nontrivial.

\brmrk
The analysis in this section shows that
if one restricts oneself to $L_2(SU(2))$,
it is not possible to get an equivariant Dirac operator
with the right summability and having a nontrivial $K$-homology class
at the same time.
Thus the classical Dirac operator for $SU(2)$,
which resides on $L_2(SU(2))\otimes\bbc^2$ is 
in some sense the minimal one.
\ermrk

\noindent{\sc Partha Sarathi Chakraborty} 
(\texttt{parthasc\_r@isical.ac.in})\\
{\footnotesize Indian Statistical 
Institute,  203, B. T. Road, Calcutta--700\,035, INDIA}\\[1ex]
{\sc Arupkumar Pal} (\texttt{arup@isid.ac.in})\\
         {\footnotesize Indian Statistical 
Institute, 7, SJSS Marg, New Delhi--110\,016, INDIA

\end{document}